\documentclass{article}
\usepackage{amsmath,amssymb,amsfonts,amsthm}
\usepackage{graphicx}

\newtheorem{theorem}{Theorem}
\numberwithin{theorem}{section}

\newtheorem{definition}[theorem]{Definition}

\newtheorem{lemma}[theorem]{Lemma}
\numberwithin{theorem}{section}

\DeclareMathOperator{\expc}{expc}
\DeclareMathOperator{\sinc}{sinc}

\begin{document}

\title{A Unified Approach to High-Gain Adaptive Controllers\thanks{This work was supported by
NSF grants \#EHS-0410685 and CMMI\#726996 as well as a Baylor University Research Council grant.}}
\author{Ian A. Gravagne\thanks{
Department of Electrical and Computer
Engineering, Baylor University, Waco, TX 76798, Email:\ Ian\_Gravagne@baylor.edu},
John M. Davis\thanks{
Department of Mathematics, Baylor University, Waco, TX 76798, Email:\
John\_M\_Davis@baylor.edu},
Jeffrey\ J. DaCunha\thanks{
Lufkin Automation, Houston, TX 77031, Email: Jeffrey\_DaCunha@yahoo.com}}
\date{}
\maketitle

\begin{abstract}
It has been known for some time that proportional output feedback will
stabilize MIMO, minimum-phase, linear time-invariant systems if the feedback
gain is sufficiently large. High-gain adaptive controllers achieve stability
by automatically driving up the feedback gain monotonically. More recently,
it was demonstrated that sample-and-hold implementations of the high-gain
adaptive controller also require adaptation of the sampling rate. In this
paper, we use recent advances in the mathematical field of dynamic equations
on time scales to unify and generalize the discrete and continuous versions
of the high-gain adaptive controller. We prove the stability of high-gain
adaptive controllers on a wide class of time scales.

\end{abstract}

\section{Introduction}

The concept of high-gain adaptive feedback arose from a desire to stabilize
certain classes of linear continuous systems without the need to explicitly
identify the unknown system parameters. This type of adaptive controller
does not identify system parameters at all, but rather adapts the feedback
gain itself in order to regulate the system. A number of papers examine the
details of various kinds of high-gain adaptive controllers \cite{Ilchmann1, Khalil, Owens2, Willems}, among others. More recently,
several papers have discussed one particularly practical angle on the
high-gain adaptive controller, namely how to cope with input/output
sampling. In particular, Owens \cite{Owens1} showed that it is not generally
possible to stabilize a linear system with adaptive high-gain feedback under
uniform sampling. Thus, Owens, et. al., develop a mechanism to adapt the
sampling rate as well as the gain, a notion subsequently improved upon by
Ilchmann and Townley \cite{Ilchmann2, Ilchmann3, Ilchmann4}, and
Logemann \cite{Logemann}.

In this paper, we employ results from the burgeoning new field of
mathematics called \textit{dynamic equations on time scales} to accomplish
three principal objectives. First, we use time scales to unify the
continuous and discrete versions of the high-gain controller, which have
previously been treated separately. Next we give an upper bound on the
system graininess to guarantee stabilizability for a much
wider class of time scales than previously known, including mixed
continuous/discrete time scales. Third, the paper represents the first
application of several very recent advances in stability theory and Lyapunov
theory for systems on time scales, and two new lemmas are presented in that
vein. We also give a simulation of a high-gain controller on a mixed time
scale.

\section{Background}

We first state two assumptions that are required in the subsequent text.

\begin{description}
\item[(A1)] The system model and feedback law are given by the linear,
time-invariant, minimum phase system
    \begin{align}
    \dot{x}(t)&=Ax(t)+Bu(t),\qquad y(t)=Cx(t),\qquad x(0)=x_{0},
    \label{system_continuous}\\
    u(t)&=-k(t)y(t),  \label{feedback_law}
    \end{align}%
for $t\ge 0$. System parameters $A\in\mathbb{R}^{n\times n}$, $B,C^{T}\in\mathbb{R}^{n\times m}$, $x_{0}\in\mathbb{R}^{n}$ and $n$ are unknown. The feedback gain $k:\mathbb{R}\rightarrow\mathbb{R}^{+}$ is piecewise continuous, and nondecreasing as $t\rightarrow \infty $. By {\it minimum phase}, we mean that the polynomial
    $$
    \det
    \begin{bmatrix}
    A-\lambda I & B \\
    C & 0
    \end{bmatrix}
    $$
with $\lambda \in\mathbb{C}$ is Hurwitz (zeros in open left hand plane).

\item[(A2)] Furthermore,
    \begin{equation}
    (CB)^{T}+(CB)>0,  \label{CB_positivedef}
    \end{equation}%
i.e. $(CB)^{T}+(CB)$ is positive definite. (In \cite{Owens2} it is pointed
out that a nonsingular input/output transformation \thinspace $T$ always
exists such that $\tilde{B}=BT^{-1}$ and $\tilde{C}=TC$ give $(\tilde{C}%
\tilde{B})^{T}+(\tilde{C}\tilde{B})>0$.)
\end{description}

Under these conditions it has been known for some time (e.g. \cite{Ilchmann1}) that there are a wide class of gain adaptation laws $k(t)=f(y(t))$, $f:\mathbb{R}^{m}\rightarrow\mathbb{R}$, that can asymptotically stabilize system (\ref{system_continuous}) in the sense that
    $$
    y\in L^{2}[t_{0},\infty ),\qquad \lim_{t\rightarrow \infty }k(t)<\infty.
    $$
Subsequently, various authors \cite{Ilchmann1, Ilchmann2, Owens1, Ozdemir} assumed that the output is obtained via sample-and-hold, i.e.
$y_{i}:=y(t_{i})$ with $k_{i}:=k(t_{i})$ and $i\in\mathbb{N}_{0}$. Thus it becomes necessary also to adapt the sample period $h_{i}:=t_{i+1}-t_{i}$ so the closed-loop control objectives are
    $$
    y_{i}\in \ell ^{2}[i_{0},\infty ),\qquad \lim_{i\rightarrow \infty
    }k_{i}<\infty ,\qquad \lim_{i\rightarrow \infty }h_{i}>0.
    $$
Though several variations on these results exist, these remain the basic
control results for continuous and discrete high-gain adaptive controllers.
The continuous and discrete cases have previously been treated quite
differently, but we now construct a common framework for both using time
scale theory.

\section{A Time Scale Model}

The system of (A1) can be replaced by%
    \begin{align}
    x^{\Delta }(t) &=\hat{A}(t)x(t)+\hat{B}(t)u(t),\qquad y(t)=Cx(t),\qquad
    x(0)=x_{0},  \label{system_timescale} \\
    u(t) &=k(t)y(t),\qquad t\in \mathbb{T},  \label{system_timescale2}
    \end{align}
where $\mathbb{T}$ is any time scale unbounded above with $0\in \mathbb{T}$.
With a series expansion similar to \cite{Ilchmann1}, we see that
    \begin{equation}
    \hat{A}(t):=\expc(\mu (t)A)A,\qquad \hat{B}(t)=\expc(\mu (t)A)B,
    \label{Ahat_Bhat}
    \end{equation}%
where expc is the matrix power series function
    \begin{equation}
    \expc(X):=I+\frac{1}{2}X+\frac{1}{6}X^{2}+\cdots+\frac{1}{n!}X^{n-1}+\cdots
    \label{expc}
    \end{equation}%
and $\mu $ is the time scale graininess.
Implementing control law \eqref{feedback_law} then gives
    \begin{equation}
    x^{\Delta }(t)=\mathcal{A}(t)x(t),\qquad \mathcal{A}(t):=\expc(\mu
    (t)A)(A-k(t)BC).  \label{system_general}
    \end{equation}
Note that $\hat{A},\hat{B}$ and $\mathcal{A}$ may all be time-varying, but
we will henceforth drop the explicit reference to $t$ for these variables.
For future reference, we also note that if $\sup_{t\in \mathbb{T}}\mu
(t)<\infty $, then $||\hat{A}||$ and $||\hat{B}||$ are also bounded (c.f. Appendix, Lemma~7.4).

The design objectives are to find graininess $\mu $ and feedback gain $k$ as
functions of the output $y$,
    \begin{align}
    \mu (t) &=g(y(t)),\qquad g:\mathbb{R}^{m}\rightarrow [0,\infty ),  \label{update_mu} \\
    k(t) &=f(y(t)),\qquad f:\mathbb{R}^{m}\rightarrow\mathbb{R}^{+},  \label{update_k}
    \end{align}
with $t\in $ $\mathbb{T}$ and $k\in C_{rd}(\mathbb{T})$ nondecreasing, such that
    $$
    y\in L^{2}[t_{0},\infty )_{\mathbb{T}}:=\left\{ w:\int_{t_{0}}^{\infty
    }|w(t)|^{2}\Delta t<\infty \right\} ,\qquad \lim_{t\rightarrow \infty
    }k(t)<\infty .
    $$

It is important to keep in mind the generality of the expressions above. A
great deal of mathematical machinery supports the existence of
delta derivatives on arbitrary times scales, as well as the existence and
characteristics of solutions to (\ref{system_general}). See, for example, e.g. \cite{Bohner1, Bohner2}.

\section{Stability Preliminaries}

We begin this section with a definition and theorem from the work of P\"{o}tzsche, Siegmund, and Wirth \cite{Potzsche}:

\begin{definition}
The {\it set of exponential stability} for the time-varying
scalar equation $z^{\Delta }(t)=\lambda (t)z(t)$, $z_{0}:=z(t_{0})\in\mathbb{C}$ with $\lambda :\mathbb{T}\rightarrow\mathbb{C}$ and $z,\lambda \in C_{rd}(\mathbb{T)}$, is given by
    $$
    \mathcal{S}(\mathbb{T}):=\mathcal{S}_{\mathbb{C}}(\mathbb{T})\cup \mathcal{S}_{\mathbb{R}}(\mathbb{T})
    $$
where
    \begin{align}
    \mathcal{S}_{\mathbb{C}}(\mathbb{T})&:=\left\{ \eta (t)\in \mathbb{C}:\alpha =-\limsup_{t\rightarrow \infty }\frac{1}{t-t_{0}}\int_{t_{0}}^{t}
\frac{\log || 1+\mu (t)\eta (t)|| }{\mu (t)}\Delta
t>0\right\},  \label{alpha_def} \\
    \mathcal{S}_{\mathbb{R}}(\mathbb{T}) &:=\left\{ \eta (t)\in \mathbb{R}:\forall t\in \mathbb{T},\exists \tau >t\text{ with }\tau \in \mathbb{T}\text{ such that }1+\mu (t)\eta (t)=0\right\},\notag
    \end{align}
with $\eta \in C_{rd}(\mathbb{T)}$ and arbitrary $t_{0}\in\mathbb{T}$.
\end{definition}

\begin{theorem} \textup{\cite{Potzsche}}
Solutions of the scalar equation $z^{\Delta }(t)=\lambda (t)z(t)$ are exponentially stable on an arbitrary $\mathbb{T}$ if and only if $\lambda (t)\in \mathcal{S}(\mathbb{T})$.
\end{theorem}

We note here that P\"{o}tzche, Siegmund, and Wirth did not explicitly consider scenarios
where $\eta $ is time-varying, but their stability analysis remains
unchanged for $\eta (t)$.

The set $\mathcal{S}_{\mathbb{R}}(\mathbb{T})$ contains nonregressive
eigenvalues $\lambda (t)$, and a loose interpretation of $\mathcal{S}_{%
\mathbb{C}}(\mathbb{T})$ suggests that it is necessary for a regressive
eigenvalue to reside in the area of the complex plane where $||
1+\mu (t)\lambda (t)|| <1$ ``most"\ of the time. The contour $%
|| 1+\mu (t)\lambda (t)|| =1$ is termed the \textit{Hilger
Circle}. Since the solution of $z^{\Delta }(t)=\lambda (t)z(t)$ is $%
x(t)=z_{0}e_{\lambda }(t,t_{0})$, Theorem 4.2 states that, if $\lambda
(t)\in \mathcal{S}(\mathbb{T})$ then some $K(t_{0})>1$ exists such that
    \begin{equation}
    || z(t)|| =|| z_{0}|| ||
    e_{\lambda }(t,t_{0})|| \leq || z_{0}||
    Ke^{-\alpha (t-t_{0})},
    \end{equation}
where $e_{\lambda }$ is a generalized time scale exponential. The Hilger Circle will be important in the upcoming
Lyapunov analysis, as will the following lemmas.

\begin{lemma}
Let $z\in C_{rd}(\mathbb{T)}$ be a function which is known
to satisfy the inequality $z^{\Delta }(t)\leq \lambda (t)z(t)$, $z(0)=z_{0}\in\mathbb{R}$ with $\lambda :\mathbb{T}\rightarrow\mathbb{R}$ and $\lambda \in C_{rd}(\mathbb{T)}$. If $z(t)>0$ for all $t\geq t_{0}\in
\mathbb{T}$, then $\lambda (t)\in \mathcal{R}^{+}$.
\end{lemma}

\begin{proof}
Defining $v(t)=z^{\Delta }(t)-\lambda (t)z(t)$ gives rise to
the initial value problem
    \begin{equation}
    z^{\Delta }(t)-\lambda (t)z(t)=v(t),\quad z(0)=z_{0},\quad t\in \mathbb{T}
    \label{y_t},
    \end{equation}
where $v(t)\leq 0$.

First, suppose $\lambda (t)$ is negatively regressive for $t\geq t_{0}$,
i.e. $\lambda (t)<-1/\mu (t)$ with $\mu (t)>0$. Then (\ref{y_t}) yields $%
z^{\sigma }=\mu \lambda z+z+\mu v=z(1+\mu \lambda )+\mu v<0$.

On the other hand, suppose $\lambda (T)$ is nonregressive for some $T>t_{0},$
$T\in \mathbb{T}$. If $\lambda (t)<-1/\mu (t)$ over $T>t\geq t_{0}$, then
invoke the preceding argument. If $\lambda (t)>-1/\mu (t)$ over $T>t\geq
t_{0}$, then solve (\ref{y_t}) \ to get
    \begin{equation}
    z(t)=e_{\lambda }(t,t_{0})z_{0}+\int_{t_{0}}^{t}e_{\lambda }(t,\sigma (\tau
    ))v(\tau )\Delta \tau ,\qquad T>t\geq t_{0}.  \label{z_regressive}
    \end{equation}
Since $e_{\lambda }(t,t_{0})>0$ for $t<T$, we see that $\int_{t_{0}}^{t}e_{%
\lambda }(t,\sigma (\tau ))v(\tau )\Delta \tau <0$. However, for $t\geq T$, (%
\ref{z_regressive}) becomes
    \begin{equation}
    z(t)=0+\int_{t_{0}}^{\rho (T)}e_{\lambda }(t,\sigma (\tau ))v(\tau )\Delta
    \tau ,\qquad t\geq T  \label{z_nonregressive}
    \end{equation}
Thus, $z(t)<0$ for all $t\geq t_{0}$ for both negatively regressive and
nonregressive $\lambda (t)$, a contradiction of the Lemma's presupposition.
This leaves only $\lambda (t)\in \mathcal{R}^{+}$.
\end{proof}

At this point we pause briefly to discuss Lyapunov theory on time scales.
DaCunha produced two pivotal works \cite{DaCunha1, DaCunha2} on solutions $P$ of the \textit{generalized time
scale Lyapunov equation},
    \begin{equation}
    A(t)^{T}P(t)+P(t)A(t)+\mu (t)A^{T}(t)P(t)A(t)=-Q(t),\text{ \ \ \ }t\in
    \mathbb{T},  \label{Generalized_Lyap}
    \end{equation}%
where $A(t),P(t),Q(t)\in\mathbb{R}^{n\times n}$, $\ A$ and $Q$ are known and $Q(t)>0$. Though it will not be
necessary to solve (\ref{Generalized_Lyap}) in this work\footnote{%
DaCunha proves that a positive definite solution $P(t)$ to the time scale
Lyapunov equation with positive definite $Q(t)$ exists if and only if the
eigenvalues of $A(t)$ are in the Hilger circle for all $t\in \mathbb{T}$.
Furthermore, $P(t)$ is unique. As with the well known result from continuous
system theory (c.f. \cite{Slotine}), the solution is constructive, with
    $$
    P(t)=\int_{{\mathbb S}_t}\Phi_{A^T(t)}(s,0)Q(t)\Phi_{A(t)}(s,0)\,\Delta s,
    $$
where $\Phi_A(t,t_0)$ denotes the transition matrix for the linear system $x^\Delta(t)=A(t)x(t)$, $x(t_0)=I$. The correct interpretation of this integral is crucial:\ for each $t\in
\mathbb{T},$ the time scale over which the integration is performed is ${\mathbb S}_t:=\mu (t){\mathbb N}_0$, which has constant graininess for each fixed $t$.}, we will see that the form of (\ref{Generalized_Lyap})
leads to an upper bound on the graininess that is generally applicable to
MIMO systems, an advancement beyond previous works which gave an explicit
bound only for SISO systems.

Before the next lemma, we define
    $$
    \mathcal{K}:=\{k:\mathbb{T\rightarrow\mathbb{R}}^{+},\ k\in C_{rd}(\mathbb{T}),\ k^{\Delta }(t)\geq 0\ \forall t\in \mathbb{T},\ \lim_{t\rightarrow \infty }k(t)=\infty \}.
    $$
The next lemma follows directly.

\begin{lemma}
Given assumptions \textup{(A1)} and \textup{(A2)} and $k\in\mathcal{K}$,
there exists a nonzero graininess $\bar{\mu}(t)$ and a time $t^{\ast }$ such
that, for all $\mu (t)\leq \bar{\mu}(t)$ and $t>t^{\ast }$, the matrix $-kC\hat{B}$ satisfies a time scale Lyapunov equation with $P=I$, $Q(t)\geq\varepsilon _{2}$ for small $\varepsilon _{2}>0$, and $\hat{B}$ from \eqref{Ahat_Bhat}.
\end{lemma}

\begin{proof}
We construct $\bar{\mu}(t)$ as
    \begin{equation}
    \bar{\mu}(t):=\frac{1}{k(t)}\left( \frac{\lambda _{\min }\{CB+(CB)^{T}\}}{\lambda _{\max }\{(CB)^{T}CB\}}-\varepsilon _{1}\right), \qquad t\in \mathbb{T},
    \label{mu_k_bound}
    \end{equation}
with $\varepsilon _{1}>0$ sufficiently small so that $\bar{\mu}(t)>0$ on $\mathbb{T}$. This holds if
    \begin{equation}
    (CB)^{T}+CB-\bar{\mu}(t)k(t)(CB)^{T}CB=\varepsilon _{1}I.
    \label{CB_Lyap2}
    \end{equation}
Multiplying (\ref{CB_Lyap2}) by $-k(t)$ gives
    $$
    (-k(t)CB)^{T}+(-k(t)CB)+\bar{\mu}(t)(-k(t)CB)^{T}(-k(t)CB)=-k(t)\varepsilon_{1}I.
    $$
(We now drop the explicit time-dependence for readability.) Set
    $$
    \Sigma(\bar{\mu}):=\expc(\bar{\mu}A)-I,
    $$
so that $\hat{B}(\bar{\mu})=[I+\Sigma (\bar{\mu})]B$, yielding
    $$
    (-kC\hat{B})^{T}+(-kC\hat{B})+\bar{\mu}(-kC\hat{B})^{T}(-kC\hat{B})=-k(\varepsilon _{1}I-Z),
    $$
where each term of $Z$ is a product of contants times $\Sigma (\bar{\mu})$.
Since $|| \Sigma (\bar{\mu})|| $ $\rightarrow 0$ as $\bar{%
\mu}\rightarrow 0$ (c.f. Appendix, Lemma 7.4), there exists a time $%
t^{\ast }\in \mathbb{T}$ when $|| Z|| <\varepsilon _{1}$.
Because the preceding arguments admit any graininess $\mu (t)\leq \bar{\mu}%
(t)$, it follows that
    \begin{equation}
    (-kC\hat{B})^{T}+(-kC\hat{B})+\mu (-kC\hat{B})^{T}(-kC\hat{B})=-kQ,\qquad
    t>t^{\ast },\ \mu (t)\leq \bar{\mu}(t),  \label{CB_Lyap3}
    \end{equation}
where $Q(t)\geq \varepsilon _{2}$ and $\varepsilon _{2}:=\varepsilon
_{1}-|| Z|| >0$.
\end{proof}

We comment briefly on the intuitive implication of Lemma 4.4. Equation \eqref{CB_Lyap3} shows that there exist positive definite $Q(t)$ and $P=I$ to
satisfy an equation of the form \eqref{Generalized_Lyap}. According to
DaCunha, this implies that the eigenvalues of $-kC\hat{B}$ lie strictly
within the Hilger circle for $t>t^{\ast }$.

\section{System Stability}

We now come to the three central theorems of the paper. If $BC$ is not known
to be full rank, or cannot be full rank because of the input/output
dimensions, then it must be assumed (or determined \textit{a priori}) that
the eigenvalues of $A-k(t)BC$ attain negative real parts at some point in
time. This phenomenon is investigated in-depth by other authors \cite%
{Ilchmann1, Willems}. We then make use of the observation by Owens, et.
al. \cite{Owens2}, that there must exist some $k^{\ast }>0$ such that if $%
k(t)\equiv k^{\ast }$ the system of (A1) has a positive-real realization.
This, together with the Kalman-Yakubovich Lemma \cite{Slotine}, implies existence of $%
P,Q>0$ so that
    \begin{equation}
    (A-k^{\ast }BC)^{T}P+P(A-k^{\ast }BC)=-Q
    \quad\text{ and }\quad
    PB=C^{T}CB.
    \label{KY}
    \end{equation}

\begin{theorem} [Exponential Stability] In addition to \textup{(A1)} and \textup{(A2)},
suppose
    \begin{enumerate}
    \item[\textup{(i)}] $t\in \mathbb{T}$ where ${\mathbb T}$ is a time scale which is unbounded above but with $\mu(t)\leq \bar{\mu}(t)$,
    \item[\textup{(ii)}] $k\in \mathcal{K}$ \textup{(}implying from \eqref{mu_k_bound} that $\mu(t)\rightarrow 0$, but not necessarily monotonically\textup{)},
    \item[\textup{(iii)}] $BC$ is not necessarily full rank, but there exists a time $t^{\ast
}\in \mathbb{T}$ such that the eigenvalues of $(A-k(t)BC)$ are strictly in
the left-hand complex plane for $t\geq t^{\ast }$.
    \end{enumerate}
Then the system \eqref{system_timescale}, \eqref{system_timescale2} is
exponentially stable in the sense that there exists time $t_{0}\in \mathbb{T}$ and constants $K(t_{0})\geq 1$, $\alpha \geq 0$, such that
    $$
    || x(t)|| \leq || x(t_{0})|| Ke^{-\frac{1}{2}\alpha
(t-t_{0})}, \qquad t\geq t_{0}.
    $$
\end{theorem}

\begin{proof}
Set $k^{\ast }:=k(t^{\ast })$. Then assumption (iii) is
the prerequisite for equation (\ref{KY}). Again, we suppress the
time-dependence of $x$, $k$, $\mu $, $\hat{A}$ and $\hat{B}$. Similarly to
Lemma 4.4, terms containing $\Sigma (\mu )$ may be added to the first
equality in \eqref{KY} to obtain
    \begin{equation}
    (\hat{A}-k^{\ast }\hat{B}C)^{T}P+P(\hat{A}-k^{\ast }\hat{B}C)\leq
    -\varepsilon _{3}\qquad t\geq t^{\ast \ast },
    \label{eps3}
    \end{equation}
for some small $\varepsilon _{3}>0$. Note $t=t^{\ast \ast }$ is the
point at which terms involving $\Sigma (\mu )$ become small enough for \eqref{eps3} to hold. Defining $Z_{1}:=C^{T}C\Sigma B-P\Sigma B$, the second
equality in \eqref{KY} gives
    $$
    P\hat{B}=C^{T}C\hat{B}+Z_{1}.
    $$

Consider the Lyapunov function $V=x(t)^{T}Px(t)$ with $P$ from \eqref{KY}.
Then, using \eqref{KY}, \eqref{eps3} and Lemma 4.4,
    \begin{align*}
    V^{\Delta }&=x^{T\Delta }Px+x^{\sigma T}Px^{\Delta }\\
    &=x^{T}[\mathcal{A}^{T}P+P\mathcal{A}+\mu \mathcal{A}^{T}P\mathcal{A}]x\\
    &=x^{T}\left[ (\hat{A}-k^{\ast }\hat{B}C)^{T}P+P(\hat{A}-k^{\ast }\hat{B}C)-(k-k^{\ast})(C^{T}\hat{B}^{T}P+P\hat{B}C)\right.\\
    &\qquad \left. +\mu (\hat{A}-k\hat{B}C)^{T}P(\hat{A}-k\hat{B}C)\right] x\\
    &\leq x^{T}[-\varepsilon _{3}-(k-k^{\ast })(C^{T}\hat{B}^{T}C^{T}C+C^{T}C
    \hat{B}C+2Z_{1}C)\\
    &\qquad +\mu \hat{A}^{T}P\hat{A}-2\mu k\hat{A}^{T}P\hat{B}C+\mu k^{2}C^{T}\hat{B}^{T}C^{T}C\hat{B}C+\mu k^{2}C^{T}\hat{B}^{T}Z_{1}C]x\\
    &=x^{T}[-\varepsilon _{3}+\mu \hat{A}^{T}P\hat{A}]x\\
    &\qquad +y^{T}[-k(C\hat{B})^{T}-kC\hat{B}+\mu k^{2}(C\hat{B})^{T}C\hat{B}+\mu k^{2}\hat{B}^{T}Z_{1}+k^{\ast}((C\hat{B})^{T}+C\hat{B})]y\\
    &\qquad +x^{T}[-(k-k^{\ast })2Z_{1}-2\mu k\hat{A}^{T}P\hat{B}]y\\
    &\leq -(\varepsilon _{3}-\mu ||\hat{A}^{T}P\hat{A}||)x^{T}x\\
    &\qquad -(k\varepsilon _{2}-\mu k||\hat{B}||
    kZ_{1}||-2k^{\ast }||C\hat{B}||)y^{T}y\\
    &\qquad +(2||kZ_{1}||+2\mu k||\hat{A}^{T}P\hat{B}||)x^{T}y.
    \end{align*}
At this point we observe the following:

\begin{itemize}
\item $||\hat{A}^{T}P\hat{A}||$ is bounded because $\hat{A%
}$ is bounded. Set $\gamma _{1}:=\sup_{t\in \mathbb{T}}||\hat{A}%
^{T}P\hat{A}||.$

\item $\mu k$ is bounded by assumption (i)\ and \eqref{mu_k_bound}.

\item $||kZ_{1}||$ is bounded because $k$ is
proportional to $\frac{1}{\mu }$ and $\frac{1}{\mu }||\Sigma (\mu
)||\rightarrow const$ as $\mu \rightarrow 0$. Set $\gamma
_{2}:=\sup_{t\in \mathbb{T}}2||kZ_{1}||.$

\item Set $\gamma _{3}:=\sup_{t\in \mathbb{T}}\mu k||\hat{B}|| \gamma _{2}$.

\item Set $\gamma _{4}:=\sup_{t\in \mathbb{T}}2k^{\ast }|| C\hat{B%
}|| .$

\item Set $\gamma _{5}:=\sup_{t\in \mathbb{T}}2\mu k|| \hat{A}%
^{T}P\hat{B}|| .$
\end{itemize}

Recalling the standard inequality $x^{T}y\leq \beta x^{T}x+\frac{1}{\beta}y^{T}y$ for any $\beta >0$, we then have
    $$
    V^{\Delta }\leq -(\varepsilon _{3}-\mu \gamma _{1}-\beta \gamma _{2}-\beta
    \gamma _{5})x^{T}x-(k\varepsilon _{2}-\gamma _{3}-\gamma _{4}-\frac{1}{\beta}\gamma _{2}-\frac{1}{\beta}\gamma _{5})y^{T}y.
    $$
By assumption (ii), there exists a time $t_{0}\geq \max \{t^{\ast
},t^{\ast \ast }\}$ such that, for sufficiently small $\beta $,
    $$
    V^{\Delta }\leq \frac{-(\varepsilon _{3}-\mu \gamma _{1}-\beta \gamma_{2}-\beta \gamma _{5})}{\gamma_{\min}(P)}V:=\eta (t)V,\qquad t\geq t_{0},
    $$
with $\eta (t)<0$. Then, by \cite[Theorem 6.1]{Bohner1}, Theorem 4.2, and
Lemma 4.3 it follows that there exists $K(t_{0})\geq 1$ so that
    $$
    || x(t)|| \leq K|| x(t_{0})|| e^{-\frac{1}{2}\alpha (t-t_{0})},\qquad t\geq t_{0},
    $$
where $\alpha $ is defined in \eqref{alpha_def}.
\end{proof}

We point out that, when $BC$ is full rank, assumption (iii) above is no
longer necessary as there always exists a $k^{\ast }$ such that the
eigenvalues of $(A-k(t)BC)$ are strictly real-negative for $k(t)\geq k^{\ast
}$. One more lemma is required before the next theorem.

\begin{lemma}
If $\mathbb{T}$ is a time scale with bounded graininess
(i.e. $\mu _{\infty }:=\sup_{t\in \mathbb{T}}\mu (t)<\infty $), then
    $$
    c_{1}\int_{t_{0}}^{\infty }e^{\alpha t}\,dt\leq \int_{t_{0}}^{\infty}e^{\alpha t}\,\Delta t\leq c_{2}\int_{t_{0}}^{\infty }e^{\alpha t}\,dt,
    $$
where $c_{1}$, $c_{2}$, $\alpha \in \mathbb{R}$ and $c_{1}$, $c_{2}>0$.
\end{lemma}

\begin{proof}
Consider the case when $\alpha >0$. The process of time
scale integration is akin to the approximation of a continuous integral via
a left-endpoint sum of (variable width)\ rectangles. If the function to be
summed is increasing (as in this case), the sum of rectangular areas will
be less than the continuous integral, meaning $c_{2}=1$, $c_{1}<1$. One
estimate of the lower bound, then, follows by simply increasing $c_{1}$
until $c_{1}e^{\alpha t}$ meets one of the rectangle right endpoints which
are given by $e^{\alpha \rho (t)}$. Thus $c_{1}e^{\alpha t}\leq e^{\alpha
\rho (t)}$, or equivalently, $c_{1}e^{\alpha \sigma (t)}\leq e^{\alpha t}$.
This in turn yields $c_{1}\leq e^{-\alpha \mu (t)}$. Therefore, the most
conservative bound is given by $c_{1}\leq e^{-\alpha \mu _{\infty }}$. The
case for $\alpha <0$ can be argued similarly, leading to the lemma's
conclusion:%
    $$
    c_{1}=
        \begin{cases}
        e^{-\alpha \mu _{\infty }}, &\alpha >0, \\
        1, & \alpha \leq 0,
        \end{cases}
    \qquad
    c_{2}=
        \begin{cases}
            1, & \alpha \geq 0, \\
            e^{-\alpha \mu _{\infty }}, & \alpha <0.
        \end{cases}
    $$
\end{proof}

We are now in a position to state the main theorem of the paper.

\begin{theorem}
In addition to \textup{(A1)} and \textup{(A2)}, assume the prototypical
update law, $k^{\Delta }(t)=|| y(t)|| ^{2}$ with $k_{0}:=k(0)>0$. Then $\lim_{t\rightarrow \infty }k(t)<\infty$ and $y\in L_{2}[t_{0},\infty )_{\mathbb{T}}$.
\end{theorem}

\begin{proof}
For the sake of contradiction, assume $k(t)\rightarrow
\infty $ as $t\rightarrow \infty $. Then $k\in \mathcal{K}$. Theorem 5.1
yields $x\in L_{\infty }[t_{0},\infty )_{\mathbb{T}}$ and therefore $y\in
L_{\infty }[t_{0},\infty )_{\mathbb{T}}$. The solution for $k(t)$ is (by
\cite[Theorem 2.77]{Bohner1}),
    $$
    k(t)=k(t_{0})+\int_{t_{0}}^{t}|| y(\tau )|| ^{2}\Delta \tau
    \leq k(t_{0})+|| C|| ^{2}\int_{t_{0}}^{t}|| x(\tau)|| ^{2}\Delta \tau ,\qquad t\geq t_{0.}
    $$
In conjunction with Lemma 5.2, this allows
    \begin{align*}
    \lim_{t\rightarrow \infty }k(t) &\leq k(t_{0})+|| C||^{2}\int_{t_{0}}^{\infty }|| x(t)|| ^{2}\,\Delta t\\
    &\leq k(t_{0})+|| C|| ^{2}|| x(t_{0})||^{2}K^{2}\int_{t_{0}}^{\infty }e^{-\alpha (t-t_{0})}\,\Delta t\\
    &=k_{0}+|| C|| ^{2}|| x(t_{0})||^{2}K^{2}\int_{0}^{\infty }e^{-\alpha t}\,\Delta t\\
    &<\infty.
    \end{align*}
This contradicts the assumption, so it must be that $k(t)<\infty $ for $t\in
\mathbb{T}$. It also immediately follows that $\int_{t_{0}}^{\infty}|| y|| ^{2}\Delta t<\infty $.
\end{proof}

It seems possible that Theorem 5.3 may be improved to show that the output
is convergent, i.e. $y(t)\rightarrow 0$ as $t\rightarrow \infty $. This is
left as an open problem.

\section{Discussion}

We remark here that there is a great amount of freedom in the choice of the
update law for $k(t)$ (c.f. \cite{Owens1}). We use the simplest choice for
convenience (as do most authors); the essential arguments remain unchanged
for other choices. There is also freedom in the choice of update for $\mu
(t) $. The expression (\ref{CB_Lyap2}), can be simplified to%
    $$
    (CB)^{T}+CB-\mu (t)k(t)(CB)^{T}CB>0,\qquad t\in \mathbb{T},
    $$
for any $\mu (t)\leq \bar{\mu}(t)$. In the SISO case, this further reduces
to the expression derived by Owens \cite{Owens1}, that $\mu kCB<2$. It
requires the graininess (which may interpreted as the system sampling step
size for sample-and-hold systems)\ to share at least an inverse relationship
with the gain, but is otherwise quite unrestrictive. Ilchmann and Townley
\cite{Ilchmann3} note that $\mu (t)=\frac{1}{k(t)\log k(t)}$ meets $\mu
kCB<2 $ after sufficient time without knowledge of $CB$. While previous
works have always constructed a monotonically decreasing step size, the time
scale-based arguments in this paper reveal even greater freedom:\ $\mu (t)$
may actually increase, jump between continuous and discrete intervals, or
even exhibit bounded randomness. Two examples of the usefulness of this
freedom come next.

For the first example, we point out that the notation in the previous
sections somewhat belies the fact that the system's time domain (its time
scale) may be fully or partially discrete, and thus there is no guarantee in
Theorem 5.3 that the output has stabilized between samples. As pointed out
in \cite{Ilchmann3} and elsewhere, a sampled system with period $\mu $ is
detectable if and only if
    \begin{equation}
    \frac{\lambda _{k}-\lambda _{l}}{2\pi j}\mu \notin \mathbb{Z}\text{ for any }%
    \lambda _{k}\neq \lambda _{l}\text{; \ }j=\sqrt{-1},
    \label{sampled_eigvals}
    \end{equation}%
where $\lambda _{k},\lambda _{l}$ $\in \{0,spec(A)\}$. We next comment on
how to circumvent the intrasample stabilization problem. Both of the
following methods essentially permit the graininess $\mu $ to ``wiggle" a bit
so that an output sample must eventually occur away from a zero-crossing.
Recalling that $\mu (t)\leq \bar{\mu}(t)$ from Theorem 5.1(i), let $\mu (t)=%
\bar{\mu}(t)v(t)$, where $0<v\leq 1\ $is one of the sequences below:

1.) Let $v$ be consist of an infinitely repeated subsequence with $n!+1$
elements that are random numbers between $0$\ and $1$. Let these elements,
labeled $v_{1},v_{2},...,$ be irrational multiples of each other. Assume $%
\mu (t)$ converges to $\bar{\mu}v_{1}>0$ such that $y(t)=0$ but the true
continuous output is nonzero. This implies that there exist integers $k,l$
such that $\frac{(\lambda _{k}-\lambda _{l})\bar{\mu}v_{1}}{2\pi j}\in
\mathbb{Z}$. As sequence $v$ advances, there may be at worst $n!$
combinations of $j,k$ such that $\frac{(\lambda _{j}-\lambda _{k})\bar{\mu}%
v_{r}}{2\pi j}\in \mathbb{Z}$ for $r=1...n!$. However, at the next instant
in time, $\frac{(\lambda _{k}-\lambda _{l})\bar{\mu}v_{n!+1}}{2\pi j}$ must
be irrational and therefore not in $\mathbb{Z}$. The controller will detect
a nonzero output and continue to adapt $k(t)$ and $\mu (t)$. In practice, of
course, it is not possible to obtain a sequence of truly irrational numbers,
but most modern computer controllers have enough accuracy to represent the
ratio of two very large integers, so that this technique would only fail for
impractically high magnitudes of $|| \lambda _{k}-\lambda
_{l}|| $.

2.)\ Let $v$ be a sequence of random numbers in the specified range. Even in
a computer with only 8-bit resolution for $v$, the probability of $\frac{%
(\lambda _{1}-\lambda _{2})\bar{\mu}v}{2\pi i}\in \mathbb{Z}$ drops
drastically after a few sample periods.

We remark here that, if \eqref{sampled_eigvals} holds then $(\mathcal{A},C)$
is detectable because $(A,C)$ is detectable. Thus, the stability of $y$
implies the stability of $x$. We do not dwell on this here, but see e.g.
\cite{Ilchmann3} for a similar argument.

For the second example, we consider a problem posed by distributed control
networks (c.f. \cite{Cervin, Gravagne}). Here, one communication
network supports many control loops as well as a certain volume of unrelated
high-priority traffic. The (unpredictable) high-priority traffic may block
the control traffic at times, forcing longer-than-anticipated sample
periods. In normal operation, the controller may sample fast enough to
behave, for all practical purposes, like a continuous control with $\mu
(t)=0 $. At the blocking instant $t_{b}$, $\mu (t_{b})$ rises to some
unpredictable level. The {\it scheduling question} is, when should the blocked
controller emit a communications packet of high enough priority to override
the block?\ One answer is straightforward:\ just before $\mu (t_{b})$
exceeds $\bar{\mu}(t_{b})$, the known maximum sample period that will
guarantee plant stability. (In \cite{Gravagne} we suggest that even longer
delays are possible under certain conditions.)\ Intuitively, lower gains
permit longer blocking delays.

\begin{figure}
\centering
\includegraphics[scale=.5]{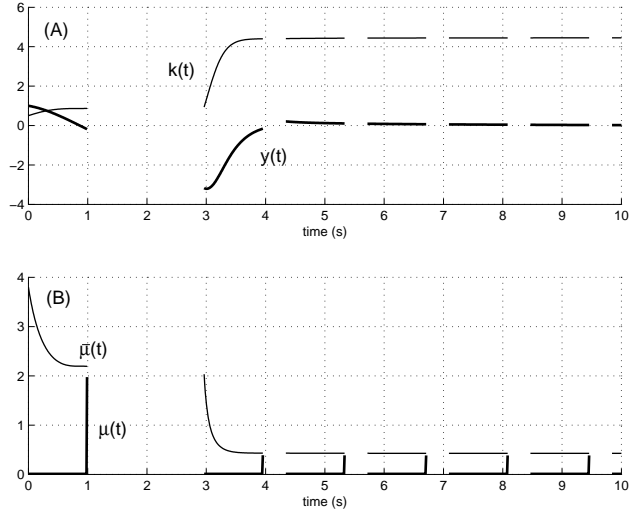}
\caption{In graph (A), output $y(t)$ (thick line) and gain $k(t)$ are plotted. In
graph (B), $\protect\mu (t)$ (thick line)\ is kept at 90\% of $\bar{\protect%
\mu}(t)$ during blocking. A block occurs after every 1 second of continuous
control.}
\label{FIG_adapt}
\end{figure}

The example of the previous paragraph is closely modeled using the variable $%
\mathbb{P}[a,b]$ time scale, which is continuous for an interval $a$, then
has a gap for an interval $b$, then repeats. Figure~\ref{FIG_adapt} shows
the regulation of a system implementing an adaptive gain controller in a
blocking situation:%
    \begin{equation}
    A=
    \begin{bmatrix}
    0 & 1 \\
    -1 & 1%
    \end{bmatrix},
    \qquad B=\begin{bmatrix}
    1 \\
    1%
    \end{bmatrix},
    \qquad C=[1,0].
    \end{equation}%
The gain begins at $k_{0}=\frac{1}{2}$ and the sampling period at $\mu (0)=0$. The bounding function for the graininess is $\bar{\mu}(t)=\frac{1.9}{k(t)}$ (so that $\bar{\mu}kCB<2$).

\section{Conclusions}

In summary, the paper illustrates a new unified continuous/discrete approach
to the high-gain adaptive controller. Using recent developments in the new
field of time scale theory, the unified results reveal that this type of
feedback control works well on a much wider variety of time scales than
explored in previous literature, including those that switch between
continuous (or nearly continuous) and discrete domains or those without
monotonically decreasing graininess. Furthermore, several results relating
to Lyapunov analysis on time scales appear here for the first time,
including Lemma 4.3 (and its use in the proof of Theorem 5.1) and 5.2. A
simulation of an adaptive controller on a mixed continuous/discrete time
scale is also given. It is our hope that time scale theory may find wider
application in the broad fields of signals and systems as it seems that many
of the tools needed in those fields are beginning to appear in their
generalized forms.

We thank our colleague, Robert J. Marks II, for his very helpful suggestions throughout this
project.

\section{Appendix}

We comment on the properties of the ``expc" function referenced in
the main body of the paper.

\begin{lemma}
The power series \eqref{expc} has the following properties:
    \begin{enumerate}
    \item $\expc(X)X=X\expc(X).$

    \item $\expc(X)=(e^{X}-I)X^{-1}$ when $X^{-1}$ exists.

    \item For real, scalar arguments $x$, $\expc(x)=e^{i\frac{x}{2}}\sinc(\frac{x}{2})$,
    where sinc denotes the sine cardinal function. \textup{(}This is the motivation for
    the expc notation.\textup{)}

    \item $|| \expc(X)|| \leq \exp (||X|| )$.

    \item $|| \expc(X)-I|| \leq \exp (||X|| )-1$.
    \end{enumerate}
    \end{lemma}

\begin{proof}
Parts 1-3 follow immediately from the definition. To verify
4, note
    $$
    \left\Vert\sum_{n=1}^{\infty }\frac{X^{n-1}}{n!}\right\Vert \leq
    \sum_{n=1}^{\infty }\frac{|| X|| ^{n-1}}{n!}\leq
    \sum_{n=1}^{\infty }\frac{|| X|| ^{n-1}}{(n-1)!}=\exp
    (|| X|| ).
    $$
Part 5 follows from a similar argument. Note that, by property 5, the
decomposition $\expc(\mu X)=I+\Sigma (\mu )$ gives $|| \Sigma
(\mu )|| \rightarrow 0$ as $\mu \rightarrow 0$, with uniform
convergence.
\end{proof}

\end{document}